\documentclass[a4paper, 11pt]{amsart}
\usepackage{amssymb,amsfonts,amsmath,amsthm,amscd,soul}

\newcommand{\cL}{\mathcal{L}}
\newcommand{\cM}{\mathcal{M}}
\newcommand{\NN}{\mathbb{N}}
\newcommand{\RR}{\mathbb{R}}
\newcommand{\CC}{\mathbb{C}}
\newcommand{\HH}{\mathbb{H}}

\newcommand{\Hom}[2]{\mathrm{Hom}(#1,#2)}
\newcommand{\End}[1]{\mathrm{End}(#1)}
\newcommand{\Ker}[1]{\mathrm{Ker}(#1)}
\renewcommand{\Im}[1]{\mathrm{Im}(#1)}
\newcommand{\Ann}[1]{\mathrm{Ann}(#1)}
\newcommand{\Ke}[1]{\mathrm{Ke}(#1)}
\newcommand{\Soc}[1]{\mathrm{Soc}(#1)}
\newcommand{\Jac}[1]{\mathrm{Jac}(#1)}
\newcommand{\Aut}{\mathrm{Aut}}
\newcommand{\DC}{\mathcal{DC}}
\theoremstyle{plain}

\newtheorem*{lem}{Lemma}
\newtheorem*{teo}{Theorem}

\newtheorem*{cor}{Corollary}
\newtheorem*{pro}{Proposition}
\newtheorem*{exa}{Example}

\title{Chain coalgebras and distributivity}
\author{Christian Lomp}
\address{Departamento de Matem\'{a}tica Pura, Universidade do Porto, Portugal}
\email{clomp@fc.up.pt}
\author{Alveri Sant'Ana}
\address{Instituto de Matem\'{a}tica, Universidade Federal do Rio Grande do Sul, Brazil}
\email{alveri@mat.ufrgs.br}

\thanks{This work was carried out as part of the project {\it Interac\c{c}\~{o}es entre \'{a}lgebras e co-\'{a}lgebras} between the Universidade do Porto and Universidade Federal do Rio Grande do Sul and Universidade de S\~{a}o Paulo financed through GRICES (Portugal) and CAPES (Brazil) during the second author's visit to the UP. 
The first author was partially supported by Centro de Matem\'{a}tica da Universidade do Porto (CMUP), financed by FCT
(Portugal) through the programs POCTI (Programa Operacional Ci\^{e}ncia, Tecnologia, Inova\c{c}\~{a}o) and POSI
(Programa Operacional Sociedade da Informa\c{c}\~{a}o), with national and European community structural funds.  The second author was supported by CAPES (Projeto 135/05 - BEX 3141/05-5) and would like to thank the Departamento de Matemática Pura for its hospitality.}

\begin{document}



\maketitle

\begin{abstract}
We show that coalgebras whose lattice of right coideals  is distributive are coproducts of coalgebras whose lattice of right coideals  is a chain. Those chain  coalgebras are characterized as finite duals of noetherian chain rings whose residue field is a finite dimensional division algebra over the base field. They also turn out to be coreflexive and infinite dimensional chain coalgebras turn out to be finite duals of left noetherian chain domains. Given any finite dimensional division algebra $D$ and $D$-bimodule structure on $D$ we construct a chain coalgebra as a cotensor coalgebra. Moreover if $D$ is separable over the base field, every chain coalgebra of type $D$ can be embedded in such a cotensor coalgebra. As a consequence cotensor coalgebras arising in this way are the only infinite dimensional chain coalgebras over perfect fields. Finite duals of power series rings with coeficients in a finite dimensional division algebra $D$ are further examples of chain coalgebras, which also can be seen as the tensor product of $D^*$ and the divided power coalgebra and can be realized as the generalized path coalgebra of a loop. If $D$ is central, any chain coalgebra is a subcoalgebra of the finite dual of $D[[x]]$. 
\end{abstract}

\section{Introduction}

Coalgebraic structures in ring theory have gained a lot of interest in recent years due to the intensive study of quantum groups and their actions and coactions. In the spirit of \cite{A} , \cite{CGT2}, \cite{FR}, \cite{JMR}, \cite{LR} and \cite{NTZ}  we are looking at the structure of coalgebras from a module theoretic point. Our aim is to describe those coalgebras whose lattice of right coideals is distributive. This is obviously the case if the lattice is linearly ordered, i.e. a chain. An important example is the  divided power coalgebra.

\begin{exa} \label{potencias divididas} 
The divided power coalgebra is the coalgebra $\DC$ with basis $\{c_n : n\in \mathbb{N} \}$ whose comultiplication and counit are given by $\Delta (c_n) = \displaystyle \sum_{i=0}^n c_i \otimes c_{n-i}$ and $\varepsilon (c_n)=\delta_{0,n}$, respectively. The coalgebra $\DC$ is cocommutative and the only subcoalgebras are given by the coradical filtration $\{\DC_n\}_{n\in\NN}$, where $\DC_n$ is the cyclic comodule generated by $c_n$ with basis $\{c_0, c_1, \ldots, c_n \}$. 
As a coalgebra $\DC_n \simeq \left( k[x]/\langle x^{n+1} \rangle\right)^*$ for any $n\geq 0$. 
It is well-known that the divided power coalgebra $\DC$ is the finite dual of the power series ring over $k$ in one variable, i.e. $\DC\simeq k[[x]]^\circ$ (see \cite{DNR}).
\end{exa}

We will call coalgebras $C$ whose lattice of right coideals is a chain, chain coalgebras. In other words, a coalgebra $C$ is a chain coalgebra if it is a uniserial left or right $C^*$-module. In \cite{CGT2} a coalgebra $C$ is called right serial if each indecomposable injective right $C$-comodule is uniserial. Moreover a right and left serial coalgebra $C$ is called uniserial if the composition factors of each indecomposable injective comodule are isomorphic. Hence chain coalgebras are uniserial coalgebras. We will show that any distributive coalgebra over a field $k$ is a coproduct of chain coalgebras and hence is serial in the sense of \cite{CGT2}. Furthermore chain coalgebras occur as finite duals of noetherian chain rings whose residue field is a finite dimensional division algebra over the base field. Duals of infinite dimensional chain coalgebras are left noetherian chain domains.

\section{Preliminaries}
Let us recall some basic notions of coalgebras: Fix a base field $k$. A (coassociative, counital) coalgebra over $k$ is a $k$-vector space  $C$ with  $k$-linear maps $\Delta:C\rightarrow C\otimes C$ (the comultiplication) and $\varepsilon:C\rightarrow k$ (the counit) such that
$$ (id\otimes \Delta)\Delta = (\Delta \otimes id)\Delta \mbox{ and } 
(id \otimes \varepsilon)\Delta = id = (\varepsilon \otimes id)\Delta$$
hold. A right $C$-comodule is a $k$-vector space $M$ with a k-linear map $\delta: M \rightarrow M\otimes C$ (the coaction) such that $(id\otimes \Delta)\delta = (\delta \otimes id)\delta$. The class of right $C$-comodules form a category $\cM^C$. Right $C$-subcomodules of $C$ are called right coideals of $C$. Moreover
attached to a coalgebra $C$ is the dual algebra $C^* = \Hom{C}{k}$ whose multiplication is the so-called {\it convolution product}, i.e. for any $f,g \in C^*$ we define $f\ast g \in C^*$ by 
$$ f\ast g  = \eta (f\otimes g) \Delta$$
where $\eta:k\otimes k \rightarrow k$ denotes the multiplication of $k$. Explicitly if $\Delta(x) = \sum_{i=1}^n x_i \otimes {x'}_i$ for some $x\in C$, then 
$$f\ast g (x)  = \sum_{i=1}^n f(x_i)g({x'}_i).$$
The unit of $C^*$ is $\varepsilon$. $C^*$ is an associative algebra. Moreover any right $C$-comodule $M$ becomes a left $C^*$-module: Let $\delta$ denote the coaction of $C$ on $M$ then for any $m\in M$ with $\delta(m)=\sum_{i=0}^n m_i \otimes c_i \in M\otimes C$ and for any $f\in C^*$ set
$$ f\rightharpoonup m = (id \otimes f)\delta(m) = \sum_{i=0}^n m_if(c_i).$$
In particular $C$ is a left $C^*$-module and the category of right $C$-comodules $\cM^C$ can be considered as a subcategory of $C^*$-Mod. More precisely $\cM^C = \sigma[_{C^*}C]$ (see \cite{BW}) where $\sigma[M]$ for a left $R$-module $M$ over some ring $R$ denotes the Wisbauer category\footnote{here we follow a recent suggestion made by Patrick F. Smith} of $M$, i.e. the subcategory of $R$-Mod whose objects are submodules of factor modules of direct sums of copies of $M$.
Analogously the category of left $C$-comodules $^C\cM$ can be considered as the subcategory $\sigma[C_{C^*}]$ of Mod-$C^*$. 
A very important fact is that $C$ is an injective cogenerator in the categories $\cM^C$ and $^C\cM$. Hence $C$ is an injective object in $\cM^C$ and any right $C$-comodule embeds into a direct sum of copies of $C$. 

Furthermore the endomorphism ring of $C$ as left $C^*$-module resp. right $C^*$-module are anti-isomorphic resp. isomorphic to $C^*$:
$$\End{_{C^*}C}^{op} \simeq C^* \simeq \End{C_{C^*}}.$$

\section{Distributive lattices}

Let $\cL= (L,\wedge, \vee, 0 ,1)$ be a complete modular lattice. Then $\cL$ is called {\it distributive} if for all elements $a,b,c \in \cL:$ $$a\wedge (b\vee c) = (a\wedge b) \vee (a\wedge c).$$
Note that the notion of distributivity is self-dual, i.e. $\cL$ is distributive if and only if $\cL^o$ is.
Attached to a left $R$-module $M$ over a unital associative ring $R$ is the lattice $\cL(M)$ of submodules of $M$. $M$ is called distributive if  $\cL(M)$ is a distributive lattice. Analogously we say that a right $C$-comodule $M$ over a coalgebra $C$ is {\it distributive} if its lattice of subcomodules is distributive. 

\subsection{}
Suppose $f: L \rightarrow L'$ is an injective homomorphism of complete modular lattices. If $L'$ is distributive, then also $L$ is.
We will use this simple argument to show that an injective cogenerator $M$ is distributive provided its endomorphism ring is. Let $M$ be a left $R$-module and $S=\End{M}$ be its endomorphism ring. Denote by
$$\mathrm{Ann}: \cL(_RM) \longrightarrow \cL(S_S)$$
the map with $\Ann{N}=\{f\in S \mid N\subseteq\Ker{f}\}$ for all submodules $N\in \cL(_RM)$. Here $\cL(_RM)$ denotes the lattice of submodules of $M$ and $\cL(S_S)$ the lattice of right ideals of $S$.

\begin{lem}\label{self-injective_homomorphism}
If $M$ is self-injective, then $\mathrm{Ann}$ is an anti-homomorphism of lattices.
\end{lem}

\begin{proof}
By definition we have $\Ann{N + L}=\Ann{N} \cap \Ann{L}$ and $\Ann{N\cap L} \supseteq \Ann{N} + \Ann{L}$ for all 
$N,L\in \cL(M)$. Let $0\neq f\in \Ann{N\cap L}$, i.e. $f(N\cap L)=0$. Then $f$ extends to a map 
\[\begin{CD}
 0 @>>> M/(N\cap L) @>>> M/N \oplus M/L \\
@. @VV\overline{f}V \\
@. M
\end{CD}\]
By self-injectivity of $M$, there exist $h:M/N\oplus M/L \rightarrow M$ such that
$\overline{f}=h$ on $M/(N\cap L)$. Writing $h=g+k$ for $g\in \Ann{N}$ and $k\in\Ann{L}$, we have $f=g+k$.
Hence $\Ann{N\cap L}=\Ann{N}+\Ann{L}$ showing that $\mathrm{Ann}$ is an lattice anti-homomorphism.  
\end{proof}

\subsection{}
Our first observation says that the endomorphism ring of a self-injective distributive module is distributive.
\begin{pro}\label{self_injective_distributive}
 Let $M$ be a self-injective left $R$-module. If $M$ is distributive as left $R$-module, then $\End{M}$ is right distributive.
\end{pro}

\begin{proof}
Suppose that $M$ is a distributive left $R$-module.
Let $I, J, K$ be right ideals of $S=\End{M}$.
Take $f \in I \cap (J+K)$, say $f = g+h$, where $f\in I$, $g\in J$ and $h\in K$. 
Then $f\in fS\cap (gS+hS)$.
Since $M$ is distributive, we have
$$ \Ker{f} + (\Ker{g}\cap \Ker{h}) = (\Ker{f}+\Ker{g}) \cap (\Ker{f}+\Ker{h}).$$
Thus applying $\mathrm{Ann}$ we get
$$\Ann{\Ker{f} + (\Ker{g}\cap \Ker{h})} = \Ann{(\Ker{f}+\Ker{g}) \cap (\Ker{f}+\Ker{h})}.$$
Since $\Ann{\Ker{f}} = fS$ by \cite[p.230]{W1} and as $\mathrm{Ann}$ is an anti-homomorphism of lattices, by \ref{self-injective_homomorphism}, we have
$$fS\cap (gS + hS) = (fS\cap gS) + (fS\cap hS).$$
Thus 
$f\in (fS\cap gS) + (fS\cap hS)  \subseteq (I\cap J) + (I\cap K)$ for all $f\in I\cap (J+K)$, i.e.
$I\cap (J+K) = I\cap (J+K)$, showing that $S$ is right distributive. 
\end{proof}

\subsection{}
For a module with some cogenerator properties, the map $\mathrm{Ann}$ becomes injective. Given a module $X$ and a module $M$, one says that $M$ cogenerates $X$ if $X$ embeds into a direct product of copies of $M$. A module that cogenerates all its factor modules is called a self-cogenerator.
\begin{lem}\label{self-cogenerator_injective}
If  $M$ is a self-cogenerator, then $\mathrm{Ann}$ is injective.
\end{lem}
\begin{proof}
Let $N,L$ be two submodules of $M$ with $\Ann{N}=\Ann{L}$. Denote by $\varphi:M/N \rightarrow M^\Lambda$ the embedding of $M/N$ into a direct product of copies of $M$. Let $\pi_\lambda:M^\Lambda\rightarrow M$ be the projection onto the $\lambda$-component for any $\lambda\in \Lambda$ and let $p_N:M\rightarrow M/N$ denote the canonical projection.
Then $\forall \lambda\in \Lambda, \pi_\lambda\varphi p_N\in\Ann{N} = \Ann{L}$ and hence for any $l\in L$ and $\lambda\in \Lambda$ we have $\pi_\lambda\varphi p_N(l)=0$. Thus $\varphi p_N(l)=0$ and as $\varphi$ was injective, $p_N(l)=0$ or equivalently $l\in N$ showing $L\subseteq N$. Analogously one shows $N\subseteq L$, i.e. $N=L$.
\end{proof}

\subsection{}
Combining the last two Lemmas we showed that the lattice of submodules of a left $R$-module which is self-injective and a self-cogenerator is isomorphic to a sublattice of the dual lattice of right ideals of its endomorphism ring. 
Hence we conclude

\begin{cor}\label{Cor_DistributiveEndoRing}
 Let $M$ be a self-injective left $R$-module which is a self-cogenerator. Then $M$ is a distributive left $R$-module if and only if $\End{M}$ is right distributive. 
\end{cor}

\begin{proof}
 By \ref{self-injective_homomorphism} and \ref{self-cogenerator_injective} if $\End{M}$ is right distributive, then $M$ is a distributive left $R$-module. The converse follows from \ref{self_injective_distributive}.
\end{proof}

\subsection{}
Considering a coalgebra $C$ as a left $C^*$-comodule and identifying the opposite ring of its endomorphism ring with $C^*$ we conclude:
\begin{cor}
The following statements are equivalent for a  coalgebra $C$ over a field. 
\begin{enumerate}
 \item[(a)] $C$ is a distributive right $C$-comodule.
\item[(b)] $C^*$ is a left distributive ring.
\item[(c)] $C$ is a distributive left $C$-comodule.
\end{enumerate}
\end{cor}

\begin{proof}
Follows from \ref{Cor_DistributiveEndoRing} (and its opposite version) and the ring isomorphisms 
$(C^*)^{op}\simeq \End{_{C^*}C}$ and $C^*\simeq \End{C_{C^*}}$.
\end{proof}

\section{Distributive coalgebras}

As seen in the last paragraph a coalgebra is a distributive right comodule if and only if it is a distributive left comodule. Hence we shall simply say that a coalgebra is distributive if it is a left or right distributive comodule. 

\subsection{}
Stephenson characterized distributive modules as those modules $M$ such that $\Hom{N/(N\cap L)}{L/(N\cap L)}=0$ for all submodules $N, L$ of $M$ (see \cite{St}). This important characterization will help us to show that any coalgebra which is distributive as left or right comodule over itself is a coproduct of coalgebras which are uniserial as comodules.
Recall that any coalgebra $C$ can be written as $$C=\bigoplus_{\beta \in I} E(S_\beta)$$ where $S_\beta$ are simple right $C$-comodules and $E(S_\beta)$ are their injective hulls in $\cM^C$. Moreover, since $C$ is a cogenerator in $\cM^C$ every simple right $C$-comodule is isomorphic to one of the $S_\beta$. Furthermore every right $C$-comodule has an essential socle, since comodules are locally finite dimensional. 

\subsection{}
We need the following module theoretic Lemma in the sequel:
\begin{lem}\label{semi_injective}
The endomorphism ring of a self-injective uniserial left $R$-module is right chain ring.
\end{lem}

\begin{proof}
Assume $M$ is self-injective and uniserial and let $f,g \in S=\End{M}$ be two endomorphisms.
 Since $M$ is uniserial we have $\Ker{f}\subseteq \Ker{g}$ or $\Ker{g}\subseteq \Ker{f}$. Suppose $\Ker{f}\subseteq \Ker{g}$, then $\Ann{\Ker{g}} \subseteq \Ann{\Ker{f}}$. By \cite[p.230]{W1} $gS = \Ann{\Ker{g}}$ and hence $gS\subseteq fS$, i.e. $S$ is a right chain ring.
\end{proof}

\subsection{}
The next Lemma tells us that direct products of chain rings are distributive.

\begin{lem}\label{Chain_ring_lemma}
 The direct product of right chain rings is a right distributive ring.
\end{lem}

\begin{proof}
 
Suppose that $R$ is a ring which is a direct product of right chain rings, say 
$R=\displaystyle\prod_{i\in I} R_i$, where the $R_i$ are right chain ring.  We first observe that $\mathcal{M}ax_r(R)$, the set of all maximal right ideals of $R$, is of the form $\mathcal{M}ax_r(R) = \{ P_j  : j\in I\}$ where for each $j\in I$ we define:
$$P_j = \{ (x_i)_{i \in I} \in R \mid x_j \in  J(R_j)\} \}$$ where $J(R_j)$ is the Jacobson radical of $R_j$. Thus, if $P_j \in \mathcal{M}ax_r(R)$, for some $j\in I$ and we take elements $x = (x_i)_{i\in I} \in R$ and 
$s = ( s_i)_{i\in I} \not\in P_j$, then we have $s_j \not\in J(R_j)$. Hence  $s_j$ is a unit element in 
$R_j$.  Set $t=(\delta_{ji})_{i\in I}$, i.e. $t_i=0$ for all $i\neq j$ and $t_j=1$ and $b=(\delta_{ji}s_j^{-1}x_j)_{i\in I}$, i.e. $b_i=0$ for all $i\neq j$ and $b_j=s_j^{-1}x_j$.
Then $xt = sb$ since 
$$(xt)_i=x_it_i=0=s_ib_i=(sb)_i$$
for all $i\neq j$ and
$$(xt)_j = x_j = s_js_j^{-1}x_j = (sb)_j.$$
Therefore, for every $P\in\mathcal{M}ax_r(R)$, the set $S_P = R\setminus P$ is a right Ore set.

Moreover, let $x = (x_i)_{i\in I}, y = (y_i)_{i \in I} \in R$. Let $j\in I$ and
$P_j\in\mathcal{M}ax_r(R)$. Without loss of generality we suppose $x_j \in y_jR_j$, since $R_j$ is a right chain ring. Now we take again $t = (\delta_{ji})_{i\in I}$, then we have $t\not\in P_j$ and $xt \in yR$. Hence it follows that $R$ is a right chain ring by \cite[Proposition 3.1]{FS}.

\end{proof}

\subsection{}
Recall that a module is called semiartinian if every factor module has an essential socle.
Moreover a family of modules $\{M_\beta\}_{\beta \in I}$ is called unrelated if $\Hom{M_\beta/N}{M_\gamma/L}=0$ for all $N\subseteq M_\beta$ and $L\subseteq M_\gamma$ for all $\beta\neq\gamma$.
Having the properties of coalgebras mentioned above in mind we prove first the following module theoretic theorem:

\begin{teo}\label{general theorem on distributive cogenerators} Let $M$ be a semiartinian, self-injective left $R$-module which is a self-cogenerator and has an indecomposable decomposition $M=\bigoplus_{\beta \in I} M_\beta$. 
Then the following conditions are equivalent:
\begin{enumerate}
 \item[(a)] $M$ is a distributive left $R$-module.
 \item[(b)]  $\{ M_\beta \}_{\beta \in I}$ is an unrelated family of uniserial modules.
\item[(c)] $\End{M}$ is isomorphic to a direct product of right chain rings.
\item[(d)] $\End{M}$ is a right distributive ring.
\end{enumerate}
\end{teo}
\begin{proof}
First note that all modules $M_\beta$ are injective in $\sigma[M]$ and have an simple essential submodule, i.e. $M_\beta$ is  the injective hull in $\sigma[M]$ of a simple $S_\beta$. 

$(a)\Rightarrow (b)$ Assume that $M$ is distributive. For any $\beta \neq \gamma$  we have $\Hom{S_\beta}{S_\gamma}=0$ by Stephenson's characterization of distributivity. Hence all simple modules $S_\beta$ are non-isomorphic.
Let $\beta \in I$ and $N\subseteq M_\beta$. Since $M$ is semiartinian, $M_\beta/N$ has an essential socle. Moreover since $M$ is self-cogenerator, $M_\beta/N$ embeds into a direct product of copies of $M$. In particular if $T$ is 
 simple submodule of $M_\beta/N$ must be isomorphic to some simple submodule of $M$, i.e. $T\simeq S_\gamma$ for some $\gamma \in I$. Let $f:T\rightarrow S_\gamma$ be such an isomorphism. By the injectivity of $M_\gamma$, $f$ can be extended to some non-zero $g:M_\beta/N \rightarrow M_\gamma$:
\[\begin{CD}
 0 @>>> T @>>> M_\beta/N \\
@. @VfVV @VVgV \\
0 @>>> S_\gamma @>>> M_\gamma
\end{CD}\]
For the projection $\pi:M_{\beta} \rightarrow M_{\beta}/N$ we then have $g\pi \in \Hom{M_\beta}{M_\gamma}$ which is $0$ if $\beta\neq \gamma$ by Stephenson's characterization. Thus $\beta=\gamma$ showing that all simple submodules of $M_\beta/N$ are isomorphic to $S_\beta$. Again by Stephenson's characterization, the socle of $M_\beta/N$ 
has to be simple. We proved that any factor module of $M_\beta$ has an essential simple socle. 
Thus $M_\beta$ is uniserial, because for any submodules $N,L$ of $M_\beta$ we have $N/(N\cap L) \cap L/(N\cap L) = 0$ in $M/(N \cap L)$ which is indecomposable. Hence $N\subseteq L$ or $L\subseteq N$.
To show that the family $\{M_\beta\}$ is unrelated assume that $f:M_\beta/N \rightarrow M_\gamma/L$ is a non-zero homomorphism for some proper submodules $N$ and $L$ of $M_\beta$ resp. $M_\gamma$. Without loss of generality  we might assume $f$ to be injective by passing to $M_\beta/N'$ where $N'/N=\Ker{f}$. Hence $$S_\beta \simeq f(\Soc{M_\beta/N})=\Soc{M_\gamma/L} \simeq S_\gamma$$  implies $\beta=\gamma$.

$(b)\Rightarrow (c)$
Since $\{ M_\beta \}$ is an unrelated family, we have
$$\Hom{M_\beta}{M} = \End{M_\beta} \oplus \Hom{M_\beta}{\bigoplus_{\gamma\neq \beta} M_\gamma} = \End{M_\beta}$$
for all $\beta \in I$. Hence
$$ \End{M} = \Hom{\bigoplus_{\beta \in I} M_\beta}{M} \simeq \prod_{\beta \in I} \Hom{M_\beta}{M} = \prod_{\beta \in I} \End{M_\beta}.$$
By Lemma \ref{semi_injective} $\End{M_\beta}$ is right uniserial.\\

$(c) \Rightarrow (d)$
Follows from \ref{Chain_ring_lemma}.\\

$(d)\Rightarrow (a)$ follows from \ref{Cor_DistributiveEndoRing}.
\end{proof}
\subsection{}
Before we apply the last theorem to the case of coalgebras we need the following Lemma:

\begin{lem} \label{decomp subcomod in direct sum of subcoalg} 
Let $C = D\oplus E$ be a coalgebra, where $D, E$ are right coideals of $C$. If $\Hom{{_{C^*}D}}{{_{C^*}E}} = 0$, then $D$ is a subcoalgebra of $C$.
\end{lem}

\begin{proof} 
Let $x\in D$. We need to show that $\Delta(x) \in D\otimes D$. For this, we can suppose that $\Delta(x) = \sum_{i=1}^n x_i\otimes y_i \in D\otimes C$ with $\{x_i\}_{1\leq i \leq n}$ linearly independent. For every $i\in \{1,2,...,n\}$ we consider the projection $\pi_i: C \rightarrow k$, where $\pi_i(v) = \lambda$, if $v = \lambda x_i$ and $\pi_i(v) = 0$ if $v\not\in kx_i$. Thus we have $\pi_i\in C^{\ast}$, for every $1\leq i \leq n$, and 
$$x\leftharpoonup \pi_j = \sum_{i=1}^n \pi_j(x_i)y_i = y_j,$$ for every $j\in \{1, 2, ..., n\}$. Now, for every $i\in \{1, 2, ..., n\}$, we define the family of maps $\psi_i: D \rightarrow C$ by $\psi_i(z)  = z\leftharpoonup \pi_i$. Then $\{\psi_i\}$ is a family of left $C^{\ast}$-linear maps. In fact, since $C$ is a $(C^{\ast}, C^{\ast})$-bimodule, we have $$\psi_i(f\rightharpoonup z) = (f\rightharpoonup z)\leftharpoonup \pi_i = f\rightharpoonup (z\leftharpoonup \pi_i) = f\rightharpoonup \psi_i(z).$$ Therefore, 
$$\psi_i \in \Hom{{_{C^*}}D}{{_{C^*}C}} = \End{{_{C^*}D}} + \Hom{{_{C^*}D}}{{_{C^*}E}} = \End{{_{C^*}D}},$$ by assumption. Hence, for all $j\in \{1, 2, ..., n\}$, we have $\Im{\psi_j}\in D$ and 
$y_j = x\leftharpoonup \pi_j = \psi_j(x) \in D$. Thus $\Delta(x) = \sum_{i=1}^n x_i\otimes y_i \in D\otimes D$ and so $D$ is a subcoalgebra of $C$, as required.
\end{proof}

\subsection{}

Recall that the coproduct of a family of coalgebras $\{C_\lambda\}_\Lambda$ is defined to be the direct sum $C=\bigoplus_\Lambda  C_\lambda$ with componentwise comultiplication and counit, i.e. 
$$\Delta i_\lambda = \Delta_\lambda \:\: \mbox{ and }\:\: \varepsilon i_\lambda = \varepsilon_\lambda,$$
where $i_\lambda$ denotes the embedding of $C_\lambda$ in $C$ and $\Delta_\lambda$ (resp. $\varepsilon_\lambda$) denotes the comultiplication (resp. counit) of $C_\lambda$.

Thus, in view of Lemma \ref{decomp subcomod in direct sum of subcoalg}, if $C$ can be written as the direct sum of an unrelated family of right comodules $C_\lambda$, then those comodules are already subcoalgebras and $C$ is the coproduct of those coalgebras. 

\subsection{}
As a corollary we have our main characterization of distributive coalgebras. Call a coalgebra $C$ right (resp. left) chain coalgebra, if it is uniserial as right (resp. left) comodule.
\begin{teo}\label{teo estrutura1}
The following conditions are equivalent for a coalgebra $C$ over a field:
\begin{enumerate}
\item[(a)] $C$ is a distributive coalgebra.
\item[(b)] $C$ is a coproduct of right chain coalgebras;
\item[(c)] $C^*$ is isomorphic to a direct product of left uniserial rings.  
\item[(d)] $C^*$ is a left distributive ring.
\end{enumerate}
\end{teo}

\begin{proof}
 This follows from \ref{general theorem on distributive cogenerators}. Since $C$ is a direct sum of an unrelated family of uniserial comodules $M_\beta$  by \ref{general theorem on distributive cogenerators}, each of this comodules is a subcoalgebra by Lemma \ref{decomp subcomod in direct sum of subcoalg}. Hence $C$ is a coproduct of chain (sub)coalgebras. 
\end{proof}

Since the indecomposable injective right coideals of $C$ are uniserial, the class of  distributive coalgebras is a subclass of right serial coalgebras introduced in \cite{CGT2}.  

\subsection{}

To any coalgebra one associates a quiver in the following way \cite{M}: If $\mathcal{S}$ is the representative set of simple coalgebras of $C$, then the quiver $\Gamma_C$ of the coalgebra $C$ has $\mathcal{S}$ as a set of vertices and there exists one arrow between $S_1$ and $S_2$, for $S_1, S_2 \in \mathcal{S}$, if and only if $S_2\wedge S_1 \not= S_1+S_2$. Also in \cite{M} it is shown that the quiver $\Gamma_C$ is isomorphic, as a directed
graph, to the Ext quiver of simple left $C$-comodules. So, we have one arrow
between $S_1$ and $S_2$ if there exists an indecomposable left $C$-subcomodule
$X$ and an exact sequence $$0\longrightarrow S_1 \longrightarrow X \longrightarrow S_2
\longrightarrow 0.$$

Let $C$ be a distributive coalgebra and suppose that $\{ S_{\beta}\}_{\beta\in I}$ is a representative set of non-isomorphic simple right $C$-comodules. Assume that there exist one arrow between $S_{\alpha}$ and $S_{\gamma}$. Then there exists an indecomposable left $C$-subcomodule and a short
exact sequence $$0\longrightarrow S_{\alpha} \longrightarrow X \longrightarrow S_{\gamma}
\longrightarrow 0.$$ 
Since simple comodules are finite dimensional, $X$ is finite dimensional. Let $E=\mathrm{Coeff}(X)$ be the coefficient subcoalgebra of $C$ associated to $X$. Since $X$ is finite dimensional also $E$ is finite dimensional, hence $E$ is a finite dimensional distributive coalgebra and $S_\alpha, X, S_\gamma$ are right $E$-comodules. By \ref{teo estrutura1}, $E$ is a finite coproduct of chain coalgebras, i.e. $E^*$ is a finite direct product of uniserial rings $E_1^*, \ldots , E_n^*$. Thus $X=X_1\times \cdots \times X_n$ is a finite direct product of $E_i^*$-modules $X_i$.
Since $X$ was indecomposable, $X=X_i$ is an $E_i^*$-module for some $i$. Since all simple $E_i^*$-modules are isomorphic, we have $S_\alpha \simeq S_\gamma$, i.e.  $\alpha=\gamma$.
We just proved:

\begin{pro}\label{quiver1}
Let $C$ be a distributive coalgebra. Then the quiver $\Gamma_C$ only
consists of isolated points or loops.
\end{pro}

The above result in fact was obtained by Cuadra and Gom\'es-Torrecillas in
\cite{CGT2} for serial coalgebras. 




\section{Chain coalgebras}

It remains to characterize right chain coalgebras to finish the classification of distributive coalgebras.
The divided power coalgebra is a typical example of chain coalgebras whose dual algebra is a power series ring, i.e. a discrete valuation ring. In this section we show that any chain coalgebra is coreflexive  and is the finite dual of a  noetherian chain ring. 

\subsection{}

Recall that the coradical filtration $\{C_n\}_{n\in\NN}$ of a coalgebra $C$ is defined to be the Loewy series of $C$ as left $C^*$-module, i.e. $C_n$ is defined as the left $C^*$-submodule that satisfies
$C_n/C_{n-1}=\Soc{C/C_{n-1}}$ for  $(n\geq 1)$, where $C_0 = \Soc{{_{C^*}C}}$.

\begin{pro}\label{BasicFactsOnUniserialCoalgebras}
Any right chain coalgebra $C$ over a field $k$ has the following properties:
\begin{enumerate}
 \item There exists a finite dimensional division algebra $D$ over $k$ with  $(C_0)^*\simeq D$ and for all $n\geq 1$: $C_n/C_{n-1} \simeq C_0$ as right and left $C$-comodule.

 \item Every proper right or left coideal of $C$ is a finite dimensional subcoalgebra of $C$ and equals $C_n$ for some $n\geq 0$. In particular $C=\bigcup_{n\in \NN} C_n$.
 \end{enumerate}
\end{pro}

\begin{proof}
(1) Since $C$ is uniserial as right $C$-comodule, $C_0$ is a simple right $C$-comodule. Since $C_0$ is a subcoalgebra, it is simple as a coalgebra. Thus $\End{_{C_0^*}{C_0}} \simeq C_0^*$ is a (finite dimensional) division ring by Schur's Lemma. Since $C$ is uniserial and every factor of $C$ has a non-zero socle, all factors $C_n/C_{n-1}$ are simple and hence isomorphic to $C_0$.

(2) First note that $\mathrm{dim}(C_n)=(n+1)\mathrm{dim}(D)$ for all $n\geq 0$. Let $N$ be a finite dimensional proper right coideal of $C$. Since $C_0\subseteq N \neq C$, there exists a maximal number $n\in \NN$ such that $C_n\subseteq N$ and $C_{n+1}\not\subseteq N$.  As $C_{n+1}/C_n$ is simple and $C_{n+1}\not\subseteq N$, we must have $N/C_n \cap C_{n+1}/C_n =0 $ which implies $N=C_n$. We recall that the members of the coradical filtration are subcoalgebras (\cite[Lemma 3.1.10]{DNR}). Note also that any element of $C$ is contained in a finite dimensional subcomodule and hence is contained in one of the $C_n$, i.e. $C=\bigcup_{n\in\NN} C_n$. Let $N$ be an infinite-dimensional right coideal of $C$. Since $C$ is a chain coalgebra and all $C_n$ are finite dimensional, $C_n\subseteq N$ for all $n\geq 0$. Hence $C=\bigcup_{n\in \NN} C_n \subseteq N \subseteq C$, i.e. $N=C$.
\end{proof}

\subsection{}\label{injectivesAreIsos}
Let $C$ be a chain coalgebra. Note first that if $C$ is infinite dimensional, then any non-zero  colinear endomorphism is surjective. To see this let $T=\End{{_{C^*}C}}$ and $0\neq f \in T$. Suppose $C$ is infinite dimensional. If $f$ would not be surjective, then there exists $n\geq 0$ such that $\Im{f}=C_n$ which is finite dimensional. Since $\Ker{f}=C_m$ for some $m\geq 0$, we have that $\mathrm{dim}(C)=\mathrm{dim}(C_n)+\mathrm{dim}(C_m)< \infty$ - a contradiction. Thus $f$ must be surjective.

In this context, we have the following result.

\begin{teo}\label{dualChainCoalgebra}
Let $A$ the dual algebra  of an infinite-dimensional chain coalgebra. Then, there exists an element $t\in A$ such
that every element of $A$ is the form $at^n$ for some $n\geq 0$, where  $a$ is  an invertible element of $A$. In particular, $A$ is a left noetherian chain domain.
\end{teo}

\begin{proof}
Let $C$ be a chain coalgebra and denote again the ring of colinear endomorphisms of $C$ by $T$. To simplify notations we denote the dual $C^*$ of $C$ by $A$. 
There exists an anti-isomorphism of algebras $\varphi:T\rightarrow A$.
Let $g\in A$ be any non-zero element and $f$ a colinear endomorphism of $C$ such that $g=\varphi(f)$.
If $f$ is injective, then by \ref{injectivesAreIsos} $f$ is an isomorphism, and hence $g$ is invertible in $A$.
Thus assuming that $g$ is not invertible in $A$, makes $f$ non-injective and as $C$ is a chain coalgebra, there exists $n\geq 0$ such that
$\Ker{f}=C_n$. Since $\Ann{\Ke{f}}=fT$, using the anti-isomorphism $\varphi$ we have that $$Ag=\varphi(\Ann{\Ke{f}}) = \varphi(\Ann{C_n}) = {C_n}^\perp = J^{n+1}$$
by \cite[Lemma 3.1.9]{DNR}, where $J=\Jac{A}$ and $I^\perp = \{ g\in A \mid g(I)=0\}$.
Hence any non-zero, non-invertible element $g\in A$ generates some finite power of $J=\Jac{A}$.
Let $I$ be any non-zero left ideal of $A$ and set $$m=\mathrm{min}\{n \mid \exists g\in I: Ag=J^n\}.$$ 
Take any $g\in I$ with $Ag = J^m$, then $I = J^m = Ag$ is cyclic. Hence $A$ is left noetherian.

\medskip
 Since any non-zero endomorphism of $C$ is surjective by \ref{injectivesAreIsos}, for all $0\neq f,g\in T$ we have $f\circ g \neq 0$. Hence  $T$ and thus $A$ are domains.

\medskip
We showed that $A$ is a local noetherian chain domain.  Since $J$ is cyclic, there exists an element $t\in A$ such that $J=At$. As $J$ is an two-sided ideal, $AtA=At$, i.e. $tA\subseteq At$.  Note that $J^n=At^n$. To see this we have $$J^n=J^{n-1}(At)=J^{n-1}t = \cdots = Jt^{n-1}=At^n.$$ 
For any element $g\in A$, $Ag=J^n$ for some $n\geq 0$. Hence $g\in At^n$, i.e. $g=ut^n$. If $u \in J$, then $u=at$ for some $a\in A$ and $g=at^{n+1}$ shows $g\in Ag=J^n \subseteq  J^{n+1}$ - a contradiction. Thus $u\not\in J$ is invertible since $A$ is local. 

\end{proof}

\subsection{}

Next we show that chain coalgebras are coreflexive. Recall that a coalgebra $C$ is coreflexive if $C$ is naturally isomorphic to $(C^*)^\circ$. By a result of Radford this is equivalent to say that all finite dimensional left $C^*$-modules are rational, i.e. are right $C$-comodules in the natural way.

\begin{pro}\label{coreflexive}
 Any infinite-dimensional chain coalgebra is coreflexive and hence isomorphic to the finite dual of a left noetherian chain domain.
\end{pro}
\begin{proof}
Since $C^*$ is noetherian, $C$ is injective in $C^*$-Mod by \cite{BW}.
As $C$ is also an essential extension of the unique simple $C^*$-module $C_0$, it is an injective cogenerator in $C^*$-Mod. Hence any left $C^*$-module embeds into a direct product of copies of $C$ as left $C^*$-module. 
 Let $M$ be a finite dimensional left $C^*$-module, then it is finitely cogenerated as left $C^*$-module and embeds into a finite direct sum of copies of $C$. Thus $M\in \sigma[{_{C^*}C}]=\cM^C$ is rational and $C$ is coreflexive by Radford's result. Since we showed that $C^*$ is a left noetherian chain domain (see \ref{dualChainCoalgebra}), we have  $C \simeq (C^*)^\circ$ as coalgebras.
\end{proof}

\subsection{}
Combining with the previous results, we see that infinite dimensional chain coalgebras are duals of left noetherian chain domains all whose left ideals are two-sided. Since commutative noetherian chain domains are discrete valuation domains (by \cite[Prop. 8.3]{Reid}), we conclude that cocommutative infinite dimensional chain coalgebras are finite duals of discrete valuation domains.

%
%

\subsection{}
Since any right (resp. left) coideal of a right chain coalgebra $C$ is equal to a term of the coradical filtration $C_n$, which is a subcoalgebra of $C$, it is also a left chain coalgebra. Hence we will simply refer to right chain coalgebras as chain coalgebras. As mentioned in the introduction, chain coalgebras form are special class of uniserial coalgebras of \cite{CGT2} and a coalgebra $C$ is a chain coalgebra if and only if it is an indecomposable right (or left) serial coalgebra.

We say that $C$ is a chain coalgebra of type $D$ for some finite dimensional division algebra $D$ over $k$ if $D\simeq (C_0)^*$. The chain coalgebras of type $k$ are precisely the pointed chain coalgebras.

The following is an easy observation from the preceding:

\begin{cor}
The  finite dimensional  chain coalgebras of type $D$ are in correspondence to the finite dimensional chain rings $A$ whose residue field is isomorphic to $D$ and  whose left ideals are the powers of $\Jac{A}$.
\end{cor}

\subsection{} 

Let $C$ be a coalgebra and $M$ a $C$-bicomodule with left coaction $\rho^l$ and right coaction $\rho^r$. The cotensor coalgebra is defined as 
$$ T_C(M) = \bigoplus_{n\geq 0} M^{\Box_C n}$$
where $M^{\Box_C n}$ for $n>0$ denotes the $n$-fold cotensor $M\Box_C \cdots \Box_C M$ and $M^{\Box_C 0}:=C$ (see \cite{N}).
The comultiplication of $T_C(M)$ is defined as

$\Delta(x)=\Delta_C(x)$ for $x\in C$, 
$\Delta(m) = \rho^l(m) + \rho^r(m)$ for $m\in M$ and 
\begin{eqnarray*}
\Delta(m_1\Box \cdots \Box m_n) = & & 
\sum (m_1)_{-1} \otimes ( (m_1)_0 \Box m_2 \cdots \Box m_n)\\
&+& \sum_{i=1}^{n-1} (m_1 \Box \cdots \Box m_i) \otimes  (m_{i+1} \Box \cdots \Box m_n)\\
&+& \sum (m_1 \Box \cdots \Box m_{n-1} \Box (m_n)_0 )\otimes (m_n)_1,
\end{eqnarray*}
where $\rho^l(m_1)=\sum (m_1)_{-1} \otimes (m_1)_0 \in C\otimes M$ and $\rho^r(m_n)=\sum (m_n)_0 \otimes (m_n)_1 \in M \otimes C$.

It has been proved in \cite[Proposition 2.3]{CGT2}, that if $Me$ is simple or zero for any primitive idempotent $e\in C_0^*$ then $T_{C_0}(M)$ is right serial. As a consequence we have:

\begin{teo}\label{cotensor_coalgebra}
 Let $D$ be any finite dimensional division algebra over $k$ and let $M$ be any $D$-bimodule with $M\simeq D$ as $k$ vector space. Then $T_{D^*}(M^*)$ is a chain coalgebra of type $D$. If $D$ is separable over $k$, then any chain coalgebra of type $D$ embeds into some cotensor coalgebra of the form $T_{D^*}(M^*)$.
\end{teo}

\begin{proof}
 Since $M\simeq D$, $M$ is a simple left  and right $D$-module. Hence $M^*$ is naturally a simple $D^*$-bicomodule and by \cite[2.3]{CGT2} $T_{D^*}(M^*)$ is right serial. By \cite[4.4]{Wood} The coradical of $T_{D^*}(M^*)$ is $D^*$ which is a simple coalgebra, we have that $T_{D^*}(M^*)$ is a uniserial comodule and hence a chain coalgebra. This proves the first claim. Now suppose that $D$ is separable over $k$, then any chain coalgebra $C$ of type $D$ has a coseparable coradical since $C_0 \simeq D^*$ is the dual of a finite dimensional separable algebra. By \cite[4.6]{Wood} $C$ embeds into $T_{C_0}(P)$ for $P=C_1/C_0$. Since the $C_0$-bicomodule $P$ is isomorphic to $C_0$ as  right and left $C$-comodule we might choose $M=P^*$ to obtain the desired $D$-bimodule where we consider $P^*$ as $C_0^*$-bimodule in the natural way and use the isomorphism $D\simeq C_0^*$ to induce a  $D$-bimodule structure on $M$. Then
$T_{C_0}(C_1/C_0) \simeq T_{D^*}(M^*)$ proves our claim.
\end{proof}

\subsection{}
The last Theorem showed that the $D$-bimodule structures on $D$ determine all chain coalgebras in case $D$ is separable over $k$. Given a finite dimensional division algebra $D$ over $k$ and ($k$-linear) automorphism $\alpha \in \mathrm{Aut}_k(D)$ we define a new $D$-bimodule structure on $D$ by 
$a \rhd x = ax$ and $x\lhd b = x\alpha(b)$ for all $a,b,x \in D$
and denote this $D$-bimodule structure by $D_\alpha$.
Suppose that $D$ is equipped with two $D$-bimodule structures and let $D_1$ resp. $D_2$ denote $D$ with those two bimodule structures. Then we say that the two $D$-bimodule structures are equivalent if there exists a $D$-bimodule isomorphism from $D_1$ to $D_2$.

\begin{lem}
 Any $D$-bimodule structure on $D$ is equivalent to the $D$-bimodule structure of $D_\alpha$ for some $\alpha\in \mathrm{Aut}_k(D)$. 
\end{lem}

\begin{proof}
Let $M=D$ be a $D$-bimodule and denote by $\rhd$ (resp. $\lhd$) the left (resp. right) $D$-action of $D$ on $M$.
Since $_DM$ is a simple left $D$-module, there exists an isomorphism of left $D$-modules
$f:M\rightarrow D$, i.e. $$\forall a,x \in D: f(a\rhd x) = a f(x).$$
Note that $f^{-1}:D \rightarrow M$ is also a left $D$-module isomorphism, where $f^{-1}(a)=a\rhd f^{-1}(1)$ for all $a\in D$.

Define $\alpha:D\rightarrow D$ by $$\alpha(a):=f(f^{-1}(1)\lhd a) \:\:\:\forall a \in D.$$
$\alpha$ is obviously an isomorphism of $k$-vector spaces and we will show, that it is an automorphism: let $a,b \in D$
first note that $$\alpha(a)\rhd f^{-1}(1) = f^{-1}(\alpha(a)) = f^{-1}(f(f^{-1}(1)\lhd a)) = f^{-1}(1)\lhd a.$$
Hence
\begin{eqnarray*}
\alpha(a)\alpha(b) &=& \alpha(a) f(f^{-1}(1)\lhd b) \\
&=& f( \alpha(a)\rhd f^{-1}(1) \lhd b)\\ 
&=& f( ( f^{-1}(1)\lhd a) \lhd b) = f(f^{-1}(1)\lhd (ab)) = \alpha(ab). 
\end{eqnarray*}
Considering $D$ with the $D$-bimodule structure of $D_\alpha$, then  $f:M\rightarrow D_\alpha$ is a $D$-bimodule isomorphism. To see this we only need to show that $f$ is right $D$-linear. Let $x,a \in D$ then
$$ f(x\lhd a) = f( f(x) \rhd f^{-1}(1) \lhd a) = f(x) f( f^{-1}(1) \lhd a ) = f(x) \alpha(a).$$

\end{proof}

Let $D$ be a finite dimensional separable division algebra over $k$. By Theorem \ref{cotensor_coalgebra}, $T_{D^*}({D_\alpha}^*)$ are the only infinite dimensional chain coalgebras of type $D$ (up to isomorphism) and  depend only on $\mathrm{Aut}_k(D)$. Let us denote for any finite dimensional division algebra $D$ over $k$ and $\alpha \in \mathrm{Aut}_k(D)$ the cotensor coalgebra $T_{D^*}({D_\alpha}^*)$ by $\DC(D,\alpha)$. At the end of the paper, we will give two concrete non-isomorphic examples of chain coalgebras arising in this way. 

Over perfect fields, any finite dimensional division algebra is separable. Hence we conclude:

\begin{cor}
 Over a perfect field, any infinite dimensional chain coalgebra is isomorphic to a cotensor coalgebra $\DC(D,\alpha)$ where $D$ is a finite dimensional division algebra and $\alpha \in \mathrm{Aut}_k(D)$.
\end{cor}

\subsection{}

Since the only ($k$-linear) $k$-bimodule structure on $k$ is the regular action, any chain coalgebra of type $k$ is isomorphic to a subcoalgebra of $\DC(k,id)=T_k(k)$ which is isomorphic to the divided power coalgebra $\DC$ (see \cite{CGT2}). Thus we have:

\begin{cor}
 The pointed chain coalgebras over $k$ are precisely the subcoalgebras of the divided power coalgebra, i.e. the dual coalgebras $\left( k[x]/\langle x^n \rangle \right)^*$ for any $n>0$ and $k[[x]]^\circ$.
\end{cor}

In particular if the field $k$ is algebraically closed then it is perfect and any chain coalgebra is pointed thus:

\begin{cor} 
Over an algebraically closed field, the subcoalgebras of the divided power coalgebra are the  only chain coalgebras
and the coproducts of subcoalgebras of the divided power coalgebra are the only distributive coalgebras.
\end{cor}

\subsection{}
Considering $D$ as a $D$-bimodule by the regular action, i.e. $D_{id}$, we will show now that the chain coalgebra $\DC(D,id)$ is isomorphic to the finite dual of the power series ring $D[[x]]$ with coeficients in $D$.

\begin{teo}
 Let $D$ be a finite dimensional division algebra over $k$. Then $\DC(D,id) \simeq D[[x]]^{o}$ as coalgebras.
\end{teo}

\begin{proof}
 Let $D$ be a division ring which is a $k$-algebra of finite dimension and
let $\{e_1, e_2, ..., e_r\}$ be a $k$-basis of $D$ with $e_1=1$ being the unit. Thus we have that $D^{\ast}$
is a coalgebra with comultiplication given by $\Delta_{D^*}(e_i^{\ast}) = \sum_{s,t}e_i^{\ast}(e_se_t)e_s^{\ast}
\otimes e_t^{\ast}$, for $1\leq i \leq r$, and counit given by $\varepsilon_{D^*}(e_i)
= e_i^{\ast}(e_1)$.

Now we consider the coalgebra $C$ which is equal to the $k$-vector space generated by  the set $$\{ e_{1n}, e_{2n}, ..., e_{rn} : n\geq 0 \}$$ with the following
 coalgebra structure 
$$\Delta(e_{kn}) = \sum_{i+j=n} \sum_{s,t} e_k^{\ast}(e_se_t)
 e_{si} \otimes e_{tj} \:\:\: 1\leq k \leq r , \mbox{ and } n\geq 0$$
 $$\varepsilon(e_{kn}) = 0, \mbox{ if } n\geq 1, \:\:\: \mbox{ and } \:\:\:
 \varepsilon(e_{k0}) = \varepsilon_{D^{\ast}}(e_k) = e_k^{\ast}(e_1)$$
 Recall that a coalgebra $C=\bigoplus_{k\geq 0} C(k)$ is called graded if $$\Delta(C(k)) \subseteq \sum_{i+j=k} C(i)\otimes C(j) \:\:\mbox{ and }\:\: \varepsilon(C(k))=0\:\:\forall k>0.$$ Let $C(k)$ be the $k$-vector space generated by $\{e_{1k}, e_{2k},..., e_{rk}\}$.
It follows by the definition of $\Delta$ that $C = \oplus_{k\geq 0} C(k)$ is
 a graded coalgebra. Also, we have that $C(0) \cong D^{\ast}$, as a coalgebra, which is the unique
 simple subcoalgebra of $C$. Let 
 $\{C_n\}_{n\geq 0}$ denote the coradical filtration of $C$. Thus  $C_0 =
 C(0)$ and $$C_1 = \Delta^{-1}(C\otimes C_0 + C_0\otimes C) = C(0)\oplus C(1).$$
 By \cite[2.2]{CMus}, $C$ is coradically graded, i.e. $C_n=\sum_{k=0}^n C(k)$.

 By the definition of $\Delta$,  $C(1)$ is a $C(0)$-bicomodule in the natural way, which allows us to consider the  cotensor coalgebra  $T_{C(0)}(C(1))$.
The canonical projection $\pi_0:C\rightarrow C(0)$ induces an $C(0)$-bicomodule map $f:C\rightarrow C(0)$ which gives rise to a coalgebra map 
$\theta: C \rightarrow T_{C(0)}(C(1))$  by the universal property
 of the cotensor coalgebra. Moreover $\theta$ is also a graded coalgebra embedding by \cite[Proposition 2.7]{CZ}.

 Now we observe that since $C(0) \cong D^{\ast}$ as a coalgebra and $C(1)
 \cong D$ as a vector space, it follows by Theorem \ref{cotensor_coalgebra} that
 $T_{C(0)}(C(1))$ is a chain coalgebra of type $D$. Moreover, since the structure
 of $D$-bimodule of $C(1)$ is the natural one, we have that $T_{C(0)}(C(1))$
 is the cotensor coalgebra $T_{D^{\ast}}(D^{\ast})=\DC(D,id)$. As a chain coalgebra has no proper infinite dimensional subcoalgebras (see \ref{BasicFactsOnUniserialCoalgebras}) the embedding $\theta$ must be an isomorphism.
 
Now we define a map $\varphi : C^{\ast} \rightarrow D[[X]]$
 by $$\varphi(f) = \sum_{n\geq 0} \left(\sum_{k=1}^r e_kf(e_{kn})\right)X^n$$
 Clearly, $\varphi$ is an isomorphism of  $k$-vector spaces. Thus, to show that $\varphi$ is an isomorphism of $k$-algebras we only need to verify that $\varphi$ preserves products. First we observe that if $e_s, e_t$ are
 two elements of the $k$-basis of $D$ as a vector space, then we have $$e_se_t
 = \sum_{k=1}^r e_k^{\ast} (e_se_t) e_k.$$ 
 
 Suppose now that $f, g\in C^{\ast}$. Then we have 
 $$\begin{array}{rcl} \varphi(f\ast g) & = & \sum_{n\geq 0}\left( \sum_{k=1}^r e_k (f\ast g)(e_{kn}\right))X^n \\
 & = & \sum_{n\geq 0} \left(\sum_{k=1}^r e_k \left[\sum_{i+j=n} \sum_{s,t=1}^r
 e_k^{\ast}(e_se_t)f(e_{si})g(e_{tj} )\right]\right)X^n \\
  & = & \sum_{n\geq 0} \left(\sum_{i+j=n} \sum_{s,t=1}^r  \left[ \sum_{k=1}^r e_k  e_k^{\ast}(e_se_t)\right]f(e_{si})g(e_{tj} ) \right)X^n \\
& = & \sum_{n\geq 0} \left(\sum_{i+j=n} \sum_{s,t=1}^r  e_s f(e_{si}) e_t  g(e_{tj} ) \right)X^n \\ 
  & = &  \left(\sum_{i\geq 0} \sum_{s=1}^r e_s f(e_{si}) X^i  \right) \left( \sum_{j\geq 0} \sum_{t=1}^r e_t g(e_{tj}) X^j\right) \\ 
 & = & \varphi(f) \varphi(g)
\end{array}$$ 
Note also that 
$$\varphi(\varepsilon) = \sum_{n\geq 0}\left( \sum_{k=1}^r e_k \varepsilon(e_{kn})\right))X^n 
= \sum_{k=1}^r e_k e_k^\ast(e_1) = e_1 = 1.$$

Hence we proved that $\DC(D,id)^\ast \simeq C^\ast \simeq D[[x]]$. Since the chain coalgebra $\DC(D,id)$ is coreflexive, we have $\DC(D,id)\simeq (\DC(D,id)^\ast)^\circ \simeq D[[x]]^\circ$.
\end{proof}

\subsection{}

There exists an alternative way to show that $D[[x]]^\circ \simeq T_{D^*}(D^*)$. Note that for a finite dimensional division algebra $D$ we have an isomorphism of $k$-algebras $D[[x]] \simeq D \otimes k[[x]]$. Moreover taking the finite duals, we get by \cite[1.5.2]{DNR}  an isomorphism $$D[[X]]^\circ \simeq D^*\otimes \DC$$ which can be shown to be a coalgebra homomorphism. Here we identify $k[[x]]^\circ$ and $\DC$. In order to apply \cite[Proposition 2.7]{CZ} to show the embedding of $C=D^*\otimes \DC$ into the cotensor coalgebra $T_{D^*}(D^*)$ we only need to show that  $C$ is coradically graded and for this one might consider the subspaces $C(k):=D^* \otimes c_k$ for $k\geq 0$, where $\{c_n \mid n\in\NN\}$ is the basis of $\DC$ as in example \ref{potencias divididas}. Then $C= \bigoplus_k C(k)$ is a graded coalgebra and one shows that the $C_0=C(0)$ and $C_1=C(0)+C(1)$, i.e. by \cite[2.2]{CMus}, $C$ is coradically graded.

\medskip

In the light of the isomorphism $\DC(D,id)\simeq D[[x]]^\circ \simeq D^*\otimes \DC$, we raise the following\\

{\bf Conjecture:}  If $D$ is a finite dimensional division algebra over $k$ and $\alpha \in \mathrm{Aut}_k(D)$ then $\DC(D,\alpha) \simeq D[[x, \alpha]]^\circ$.\\

Here $D[[x, \alpha]]$ denotes the skew power series ring, i.e. subject to the relation $Xa = \alpha(a)X$ for all $a\in D$.

\subsection{}
For chain coalgebras of central type, of type $D$ where the center of 
$D$ is $k$, we have that all $D$-bimodule structures on $D$ are equivalent, because $D\otimes D^{op}\simeq \End{D}$ is simple artinian and hence there exists just one simple $D$-bimodule up to isomorphism.
Thus we proved the

\begin{cor}
 Let $D$ be a finite dimensional central division algebra over $k$. Then any chain coalgebra of type $D$ is isomorphic to a subcoalgebra of $D[[x]]^\circ$, the finite dual of the power series ring with coeficients in $D$.
\end{cor}

\subsection{}

Having in mind that the divided power coalgebra $\DC=k[[x]]^\circ$ is the path coalgebra of a loop, one might ask whether the chain coalgebras $D[[x]]^\circ$ also have an interpretation in terms of path coalgebras.
In \cite{LL}, generalized path coalgebras were introduced and it actually turns out that those chain coalgberas $D[[x]]^\circ$ can be considered as a generalized path coalgebra of a loop, where the coalgebra $D^*$ is attached to the vertex of the loop. We will recall some basic definitions of \cite{LL}.
Let $Q =(Q_0,Q_1,s,e)$ be a quiver, where we are denoting by $Q_0$ the set
of vertices, by $Q_1$ the set of arrows, by $s(\alpha)$ the
start vertex of $\alpha$ and by $e(\alpha)$ the end vertex of $\alpha$, for
$\alpha \in Q_1$. A path in $Q$ will be denoted by $(a|\alpha_1\alpha_2...\alpha_n|b)$,
where $\alpha_i \in Q_1$ for every $i=1,2,...,n$, and $a=s(\alpha_1), \,
e(\alpha_i) = s(\alpha_{i+1}), \, 1\leq i \leq n-1$, and $e(\alpha_n) = b$.
The length of a path is the number of arrows in it. 

Given a family of coalgebras $\mathcal{C} = \{ (C_i, \Delta_i, \varepsilon_i) | i\in Q_0\}$ of coalgebras associated to each of the vertices of $Q$, the generalized
path coalgebra is constructed in the following way. The elements of $\cup_{i\in Q_0}C_i$
will be called the $\mathcal{C}$-path of length zero, and for each $n\geq
1$, a path of length $n$ will be given by $a_1\beta_1a_2\beta_2...a_n\beta_na_{n+1}$,
where $(s(\beta_1)|\beta_1\beta_2...\beta_n|e(\beta_n))$ is a path in $Q$
of length $n$ and for each $i = 1, 2, ..., n$ $a_i\in C_{s(\beta_i)}$ and
$a_{n+1}\in C_{e(\beta_n)}$.

The generalized path coalgebra  $k(Q,\mathcal{C})$ is  defined as the quotient of
the $k$-linear space generated by set of all $\mathcal{C}$-paths of $Q$ by
the space generated by all the elements of type $$a\beta_1...\beta_{j-1}\left(\sum_{l=1}^m
\lambda_la_{jl}\right)\beta_ja_{j+1}...a_n\beta_na_{n+1} - \sum_{l=1}^m \lambda_la_1\beta_1...\beta_{j-1}a_{jl}\beta_ja_{j+1}...a_n\beta_na_{n+1}$$
where $(s(\beta_1)|\beta_1...\beta_n|e(\beta_n))$ is a path in $Q$ of length
$n$, and for each $i=1,2...,n$, $a_i\in C_{s(\beta_i)}$, $a_{n+1}\in C_{e(\beta_n)}$,
and $\lambda_l \in k$, $a_{jl}\in C_{s(\beta_j)}$ for $l = 1, 2,..., m$.,
endowed with the following coalgebra structure: For a $\mathcal{C}$-path
of length $n$ $a_1\beta_1...a_n\beta_na_{n+1}$ we define $$\Delta (a_1\beta_1...a_n\beta_na_{n+1})
= \sum_{i=1}^{n+1} \sum_{a_i}( a_1\beta_1...a_{i-1}\beta_{i-1}a_i^{\prime})\otimes (a_i^{\prime\prime}\beta_ia_{i+1}...a_n\beta_na_{n+1})$$

where $\Delta_i(a_i) = \sum_{a_i} a_i^\prime\otimes a_i^{\prime\prime}$ denotes the comultiplication of $a_i\in C_i$.
The counit is defined as $$\varepsilon(p) = \left\{\begin{array}{rl} 0, & \mbox{ if the length of
$p$ is $n>0$} \\ \varepsilon_i(p), & \mbox{ if $p\in C_i$ for some $i\in Q_0$
}\end{array}\right.$$ 


Li and Liu  showed in \cite[Proposition 2.2]{LL} that
$$k(Q,\mathcal{C}) \cong T_{k(Q_0,\mathcal{C})}(k(Q_1,\mathcal{C}))$$

Let $Q$ be a loop, i.e. $Q_0=\{ 1 \}$, $Q_1=\{ \alpha \}$, $s(\alpha)=e(\alpha)=1$, and let $\mathcal{C}= \{ D^* \}$ for a finite dimensional division algebra $D$ over $k$.
Since $Q_0$ and $Q_1$ just consist of one element, $k(Q_0, \{ D^* \}) \simeq D^*$ as coalgebras and $k(Q_1, \{ D^* \}) \simeq D^*$ as $D^*$-bicomodules.

Then $$k(Q,\{ D^*\}) \simeq T_{D^*}(D^*) \simeq \DC(D,id) \simeq D[[x]]^\circ.$$

\section{Examples of chain coalgebras}

The question arises how to construct explicitly examples of chain coalgebras that are not subcoalgebras of the divided power coalgebra.

\subsection{}\label{example1}
Choose $k=\RR$ and $D=\CC$.
Let $C$ be the real vector space with basis $\{x_n, y_n \mid n\geq 0\}$ and define the following coalgebraic structure on $C$:
$$ \Delta(x_n) = \sum_{i+j=n} x_i \otimes x_j - y_i \otimes y_j $$
$$ \Delta(y_n) = \sum_{i+j=n} x_i \otimes y_j + y_i \otimes x_j $$
and with counit $\varepsilon(x_n)=\delta_{0,n}$ and $\varepsilon(y_n)=0$ for all $n\geq 0$.
Then $C^* \simeq \CC[[z]]$ by the map 
$$f\mapsto \sum_{n=0}^\infty \left( f(x_n)+ \imath f(y_n) \right) z^n.$$
Hence $C$ is the finite dual coalgebra of $\CC[[z]]$.
\subsection{}
We will construct a family of non-cocommutative chain coalgebras.
Let $D$ be any finite dimensional division algebra over $k$ and $\varphi, \sigma \in \Aut_k(D)$ two $k$-linear automorphisms of $D$. Let $V=D$ be the $(D,D)$-bimodule whose structure is given by $a\cdot x = \varphi(a)x$ and 
$x\cdot a = x\sigma(a)$, for all $a\in D$, $x\in V$. Let $A=D\propto V$ be the trivial extension of $D$ by $V$, that is as a $k$-vector space $A=D\times V$ with multiplication
$$(a,x)\cdot (b,y) = (ab, a\cdot y + x\cdot b) = (ab, \varphi(a)y + x\sigma(b)).$$ Then $A$ is a $k$-algebra of dimension $2\mathrm{dim}(D)$. Note, that $A$ has only three left (resp. right) ideals $A$, $\Jac{A}$ and $0$, where the Jacobson radical of $A$ is square zero and equals $\Jac{A}=(0,V)$.
Also note that $A$ is commutative if and only if $D$ is commutative and $\varphi=\sigma$.
Let $C=\Hom{A}{k}=A^*$ be the dual $k$-coalgebra associated to $A$. Then $$C_0 = (\Jac{A})^\perp \simeq (A/\Jac{A})^* \simeq D^*$$ is simple and $0, C_0$ and $C$ are the only three left (resp. right) coideals  of $C$, i.e. $C$ is a chain coalgebra of type $D$ which is cocommutative if and only if $D$ is commutative and $\varphi=\sigma$.

\subsection{}
To give an explicit instance of the example above let $k=\RR$ be the real numbers and $D=\CC$ be the complex numbers.
Let $\varphi=\overline{\phantom{X}}$ be the complex conjugation and $\sigma=id$. Then $A=\CC\propto\CC$ has multilplication
$$(a,x)\cdot(b,y) = (ab,\overline{a}y+xb).$$
$A$ can be seen as the factor ring $\CC[x,\alpha]/\langle x^2\rangle$ of the Ore extension  $\CC[x,\alpha]$, where $\alpha$ denotes the complex conjugation. Hence complex scalars commute with $x$ by the relation $wx=x\overline{w}$.

Let $C=A^*$ be the dual coalgebra of $A$. Then $C$ is $4$-dimensional with basis $e,f,g,h$ corresponding to the $4$ basis elements $(1,0),(\imath,0),(0,1),(0,\imath)$ of $A$ with comultiplication

\begin{eqnarray*}
\Delta(e) &=& e\otimes e - f\otimes f\\
\Delta(f) &=& e\otimes f + f\otimes e\\
\Delta(g) &=& e\otimes g + g\otimes e + f\otimes h - h\otimes f \\
\Delta(h) &=& e\otimes h + h\otimes e - f\otimes g + g\otimes f
\end{eqnarray*}

and counit $\varepsilon(e)=1$ and $\varepsilon(f)=\varepsilon(g)=\varepsilon(h)=0$.
The $4$-dimensional chain coalgebra $C$ of type $\CC$ is not cocommutative and hence not isomorphic to a subcoalgebra of the divided power coalgebra.

\subsection{}\label{Example2}
Extending the above construction we get an infinite dimensional non-cocommutative chain coalgebras:
 Choose $k=\RR$ and $D=\CC$.
Let $C$ be the real vector space with basis $\{x_n, y_n \mid n\geq 0\}$ and define the following coalgebraic structure on $C$:
$$ \Delta(x_n) = \sum_{i+j=n} x_i \otimes x_j - (-1)^i y_i \otimes y_j $$
$$ \Delta(y_n) = \sum_{i+j=n} (-1)^i x_i \otimes y_j + y_i \otimes x_j $$
and with counit $\varepsilon(x_n)=\delta_{0,n}$ and $\varepsilon(y_n)=0$ for all $n\geq 0$.
Let $\alpha$ denote the complex conjugation and let $\CC[[z, \alpha]]$ be the skew power series ring where $zw=\alpha(w)z$
for any $w\in\CC$.
Then $C^* \simeq \CC[[z, \alpha]]$ by the map 
$$f\mapsto \sum_{n=0}^\infty \left( f(x_n)+ \imath f(y_n) \right) z^n.$$
Hence $C$ is the finite dual coalgebra of the (non-commutative) discrete valuation domain  $\CC[[z,\alpha]]$.

\subsection{}
Since the only finite dimensional division algebras over $\RR$ are $\RR, \CC$ and the quaternions $\HH$ and since $\mathrm{Aut}_\RR(\CC)= \{id ,\alpha\}$ and $\HH$ is central, we have

\begin{cor}
  $\RR[[x]]^\circ, \CC[[x]]^\circ, \CC[[x, \alpha]]^\circ$ and $\HH[[x]]^\circ$ are the only real infinite dimensional chain coalgebras.
\end{cor}

\end{document}